\documentclass[11pt]{amsart}
\usepackage{amsmath}
\usepackage{amssymb}
\usepackage{amscd}
\usepackage{hhline}
\usepackage{amsfonts}%

\newtheorem{theorema}{Theorem}
\newtheorem{dfn}[theorema]{Definition}

\newtheorem{lem}[theorema]{Lemma}

\newcommand\C{{\mathbb C}}
\renewcommand\d{\delta}
\newcommand\E{{\mathcal E}}

\newcommand\La{\Lambda}
\newcommand\op[1]{\mathop{\rm #1}\nolimits}
\newcommand\ot{\otimes}

\newcommand\R{{\mathbb R}}

\newcommand\z{\sigma}
\newcommand\Z{{\mathbb Z}}

\begin{document}

 \title{Dimension of the solutions space of PDEs}
 \author{Boris Kruglikov, Valentin Lychagin}
 \address{Institute of Mathematics and Statistics, University of Troms\o,
 Troms\o, 90-37, Norway.}
 \email{kruglikov@math.uit.no, \quad lychagin@math.uit.no.}
 \subjclass[2000]{Primary: 35N10, 58A20, 58H10; Secondary: 35A30.}
 \keywords{Solutions space, Cartan's test, Cohen-Macaulay module, involutive
system, compatibility, formal integrability.}
 \maketitle

 \vspace{-14.5pt}
 \begin{abstract}
We discuss the dimensional characterization of the solutions space
of a formally integrable system of partial differential equations
and provide certain formulas for calculations of these dimensional
quantities.
 \end{abstract}

\section{Introduction: what is the solutions space?}

Let $\E$ be a system of partial differential equations (PDEs). We
would like to discuss the dimensional characterization of its
solutions space.

However it is not agreed upon what should be called a solution. We
can choose between global or local and even formal solutions or
jet-solutions to a certain order. Hyperbolic systems hint us about
shock waves as multiple-valued solutions and elliptic PDEs suggest
generalized functions or sections.

A choice of category, i.e. finitely differentiable $C^k$, smooth
$C^\infty$ or analytic $C^\omega$ together with many others, plays
a crucial role. For instance there are systems of PDEs that have
solutions in one category, but lacks them in another (we can name
the famous Lewy's example of a formally integrable PDE without
smooth or analytic solutions, \cite{L}).

In this paper we restrict to local or even formal solutions. The
reason is lack of reasonable existence and uniqueness theorems (in
the case of global solutions even for ODEs). In addition this
helps to overcome difficulties with blow-ups and multi-values.

If the category is analytic, then Cartan-K\"ahler theorem
\cite{Ka} guarantees local solutions of formally integrable
equations \cite{Go} and even predicts their quantity. We then
measure it by certain dimension characteristics.

If the category is smooth, formal integrability yields existence
of solutions only if coupled with certain additional conditions
(see for instance \cite{Ho}). Thus it is easier in this case to
turn to formal solutions, which in regular situations give the
same dimension characteristics. With this vague idea let us call
the space of solutions $\op{Sol}(\E)$.

With this approach it is easy to impose a topology on the
solutions space. However we shall encounter the situations, when
the topological structure is non-uniform.

To illustrate the above discussion, let's consider some model ODEs
(in which case we possess existence and uniqueness theorem). The
space of local solutions for the ODE $y'=y^2$ is clearly
one-dimensional, but the space of global solutions (continuous
pieces until the blow-up) has two disconnected continuous pieces
(solutions $y=(a-x)^{-1}$ for $a<0$ or for $a>0$) and a singular
point (solution $y=0$). Another example is the equation
$y'^2+y^2=1$, the solutions on$(-\varepsilon,\varepsilon)$ of
which form $S^1$, but the space of global solutions is $\R^1$
(both united with two singular points in $\op{Sol}(\E)$).

We would like to observe the "{}biggest"\ piece of the space of
$\op{Sol}(\E)$, so that in our dimensional count we ignore
isolated and special solutions or their families and take those of
connected components, that have more parameters in.

It will be precisely the number of parameters, on which a general
solution depends, that we count as a dimensional characteristic.
Let us discuss the general idea how to count it and then give more
specified definitions.

Note that in this paper we consider only (over)determined systems
of PDEs. Most results will work for underdetermined systems, but
we are not concerned with them.

\section{Understanding dimension of the solutions space}

Let us treat at first the case of linear PDEs systems (the method
can be transferred to non-linear case). We consider formal
solutions and thus assume the system of PDEs ${\mathcal{E}}$ is
formally integrable. We also assume the system $\E=\E_k$ is of
pure order $k$, which shall be generalized later.

Thus for some vector bundle $\pi:E(\pi)\rightarrow M$ we identify
$\E$ as a subbundle ${\E_k}\subset J^{k}(\pi)$ (see
\cite{S,Go,KLV}) and let $\E_l\subset J^l(\pi)$ be its $(l-k)$-th
prolongations, $l\ge k$. Then the fibres $\E_x^\infty\subset
J_{x}^{\infty}(\pi)$ at points $x\in M$ can be viewed as spaces of
formal solutions of $\E$ at $x\in M.$ To estimate size of
$\E^\infty_x$ we consider the spaces of linear functions on
${\mathcal{E}}_{l,x}$, i.e. the space $\E_{l,x}^*$. The
projections
$\pi_{l,l-1}:{\mathcal{E}}_{l,x}\rightarrow{\mathcal{E}}_{l-1,x}$
induce embeddings
$\pi^*_{l,l-1}:{\mathcal{E}}_{l-1,x}^*\hookrightarrow{\mathcal{E}}_{l,x}^*$,
and we have the projective limit
 $$
{\mathcal{E}}_{x}^{\ast}=\cup_{l}{\mathcal{E}}_{l,x}^{\ast}.
 $$

Remark that $\E^*$ is the module over all scalar valued
differential operators on $\pi$, while the kernel of the natural
projection $J_x^\infty(\pi) ^*\to \E_x^*$ can be viewed as the
space of scalar valued differential operators on $\pi$ vanishing
on the solutions of the PDEs system $\E$ at the point $x\in M$.
Thus elements of $\E_x^*$ are linear functions on the formal
solutions $\E_x^\infty$.

We would like to choose "{}coordinates"\ among them, which will
estimate dimension of the formal solution space. To do this we
consider the graded module associated with filtred module
$\E_x^*$:
 $$
g^*(x)=\bigoplus_{l\ge0} g_l^*(x),
 $$
where $g_l(x)$ are the symbols of the equation at $x\in M$:
 $$
g_l(x)=\E_{l,x}/\E_{l-1,x}\subset S^lT_x^*\otimes\pi_{x}
 $$
(we let $\E_l=J^l(\pi)$ for $l<k$), and reduce analysis of
${\mathcal{E}}_{x}^{\ast}$ to investigation of the symbolic module
$g_{x}^{\ast}$.

This $g^*$ is the module over the symmetric algebra $ST_xM=\oplus
S^i(T_xM)$ and its support $\op{Char}_x^{\mathbb{C}}(\E)\subset
\mathbb{P}^{\mathbb{C}}T^*_xM$ is a complex projective variety
consisting of complex characteristic vectors. The values
$\mathcal{K}_{p}$ of the symbolic module $g_{x}^*$ at
characteristic covectors $p\in {}^\mathbb{C}T_x^*\setminus0$ form
a family of vector spaces over $\op{Char}_x^{\mathbb{C}}(\E)$,
which we call \emph{characteristic sheaf}.

By the Noether normalization lemma (\cite{E}) there is a subspace
$U\subset T_xM$ such that the homogeneous coordinate ring
$ST_xM/\op{Ann}g^*(x)$ of $\op{Char}_x^{\mathbb{C}}(\E)$ is a
finitely generated module over $SU$. It follows that $g^*(x)$ is a
finitely generated module over $SU$ too.

If $g^*(x)$ is a Cohen-Macaulay module (see \cite{E}, but we
recall the definition later in a more general situation, then
$g^*(x)$ is a free $SU$-module (we called the respective PDEs
systems $\E$ Cohen-Macaulay in \cite{KL$_2$} and discussed their
corresponding reduction).

Let $\sigma$ be the rank of this module, and $p=\dim U$. By the
above discussion these numbers can be naturally called
\emph{formal functional rank} and \emph{formal functional
dimension} of the solutions space $\E^\infty_x$ at the point $x\in
M$, because they describe on how many functions of how many
variables a general jet-solution formally depends (we shall omit
the word "{}formally"\ later), or how many "{}coordinates"\ from
$\E^*_x$ should be fixed to get a formal solution.

If  the symbolic module is not Cohen-Macaulay, the module $g^*(x)$
over $SU$ is not free, but finitely generated and supported on
$\mathbb{P}^{\mathbb{C}}U^*$. Let $\mathbb{F}(U)$ be the field of
homogeneous functions $P/Q,$ where $P,Q\in$ $SU$, $Q\ne0$,
considered as polynomials on $U^*$. Thus $\mathbb{F}(U)$ is the
field of meromorphic (rational) functions on $U^*$.

Consider $\mathbb{F}(U)\otimes g^*(x)$ as a vector space over
$\mathbb{F}(U)$. Keeping the same definition for $\sigma$, let us
call the dimension of this vector space $p$ \emph{formal rank\/}
of $\E$ at the point $x\in M$.

It is clear that for Cohen-Macaulay systems the two notions
coincide. However since $g^*(x)$ over $SU$ is not free, we would
like to give more numbers to characterize the symbolic module.

Let us choose a base $e_1,\dots,e_r$ of $\mathbb{F}(U)\otimes
g^*(x)$ such that $e_1,\dots,e_r$ are homogeneous elements of
$g^*(x)$ and denote by $\Gamma_1\subset g^*(x)$ the $SU$-submodule
generated by this base. It is easy to check that $\Gamma_{1}$ is a
free $SU$-module. For the quotient module $M_1=g^*(x)/\Gamma_1$ we
have the following property:
 $$
\op{Ann}h\neq0 \text{ in }SU,\ \text{for any }h\in M_1.
 $$
Therefore $\op{Ann}M_1\ne 0$ and the support $\Xi_1$ of $M_1$ is a
proper projective variety in $\mathbb{P}^{\mathbb{C}}U^*$.

We apply the Noether normalization lemma to $\Xi_1$, we get a
subspace $U_1\subset U$, such that $M_1$ is a finitely generated
module over $SU_1$. Its rank will be the next number $p_1$ and we
also get $\sigma_1=\dim U_1$, which we can call the next formal
rank and formal dimension.

Applying this procedure several more times we get a sequence of
varieties $\Xi_i$ and numbers $(p_i,\sigma_i)$, which depends, in
general, on the choice of the flag $U\supset U_1\supset
U_2\supset\dots$ and the submodules $\Gamma_i$ of $SU_{i-1}$.

Thus we resolve our symbolic module via the exact 3-sequences
 \begin{multline*}
\hphantom{AAAA}0\to\Gamma_1\to g^*\to M_1\to0\ \text{ over }SU,\\
0\to\Gamma_2\to M_1\to M_2\to0\ \text{ over }SU_1,\qquad
\dots\hphantom{AAAA}
 \end{multline*}
(with
$\op{Supp}M_i=\op{Supp}\Gamma_{i+1}\supsetneqq\op{Supp}M_{i+1}$)
\text{etc.}
\section{Cartan numbers}\label{Cn}

In Cartan's study of PDEs systems $\E$ (basically viewed as
exterior differential systems in this approach) he constructed a
sequence of numbers $s_i$, which are basic for his involutivity
test. These numbers depend on the flag of subspaces one chooses
for investigation of the system and so have no invariant meaning.

The classical formulation is that a general solution depends on
$s_p$ functions of $p$ variables, $s_{p-1}$ functions of $(p-1)$
variables, \dots, $s_1$ functions of 1 variable and $s_0$
constants (we adopt here the notations from \cite{BCG$^3$}; in
Cartan's notations \cite{C} we should rather write $s_p$,
$s_p+s_{p-1}$, $s_p+s_{p-1}+s_{p-2}$ etc). However as Cartan
notices just after the formulation \cite{C}, this statement has
only a calculational meaning.

Nevertheless two numbers are absolute invariants and play an
important role. These are Cartan genre, i.e. the maximal number
$p$ such $s_p\ne0$, but $s_{p+1}=0$, and Cartan integer
$\sigma=s_p$.

As a result of Cartan's test a general solution depends on
$\sigma$ functions of $p$ variables (and some number of functions
of lower number of variables, but this number can vary depending
on a way we parametrize the solutions). Here general solution is a
local analytic solution obtained as a result of application of
Cartan-K\"ahler (or Cauchy-Kovalevskaya) theorem and thus being
parametrized by the Cauchy data.

Hence we can think of $p$ as of {\em functional dimension} and of
$\sigma$ as of {\em functional rank} of the solutions space
$\op{Sol}(\E)$. In fact, we adopt this terminology further on in
the paper, because as was shown in the previous section it
correctly reflects the situation.

These numbers can be computed via the characteristic variety. If
the characteristic sheaf over $\op{Char}^\C(\E)$ has fibers of
dimension $k$, then
 $$
p=\dim\op{Char}^\C(\E)+1,\ \sigma=k\cdot\op{deg}\op{Char}^\C(\E).
 $$
The first formula is a part of Hilbert-Serre theorem (\cite{H}),
while the second is more complicated. Actually Cartan integer
$\sigma$ was calculated in \cite{BCG$^3$} in general situation and
the formula is as follows.

Let $\op{Char}^\C(g)=\cup_\epsilon\Sigma_\epsilon$ be the
decomposition of the characteristic variety into irreducible
components and $d_\epsilon=\dim \mathcal{K}_x$ for a generic point
$x\in\Sigma_\epsilon$. Then
 $$
\sigma=\sum d_\epsilon\cdot\op{deg}\Sigma_\epsilon.
 $$

It is important that these numbers coincide with the functional
dimension and rank of the previous section. Moreover the sequence
of Cartan numbers ${s_i}$ is related to the sequence
${(p_i,\sigma_i)}$ of the previous section.

This can be seen from the general approach of the next and
following sections, which treat the case of systems $\E$ of PDEs
of different orders (we though make presentation for the symbolic
systems, with interpretation for general systems being well-known
\cite{S,KLV,KL$_2$}).

\section{Symbolic systems}

Consider a vector space $T$ of dimension $n$ (tangent space to the
set of independent variables, substitute to $T_xM$) and a vector
space $N$ of dimension $m$ (tangent space to the set of dependent
variables, substitute to $\pi_x=\pi^{-1}(x)$).

Spencer $\d$-complex is de Rham complex of polynomial
$N$-valued differential forms on $T$:
 $$
0\to S^kT^*\ot N\stackrel{\d}\to S^{k-1}T^*\ot N\ot
T^*\stackrel{\d}\to\dots\stackrel{\d}\to S^{k-n}T^*\ot
N\ot\La^nT^*\to0,
 $$
where $S^iT^*=0$ for $i<0$. Denote by
 $$
\d_v=i_v\circ\d :S^{k+1}T^*\ot N\to S^kT^*\ot N
 $$
the differentiation along the vector $v\in T$.

The $l$-th prolongation of a subspace $h\subset S^kT^*\ot N$ is
 \begin{multline*}
h^{(l)}=\{p\in S^{k+l}T^*\ot N\,:\,\d_{v_1}\dots\d_{v_l}p\in h\
\forall v_1,\dots,v_l\} =S^lT^*\ot h\cap S^{k+l}T^*\ot N.
 \end{multline*}

 \begin{dfn}
A sequence of subspaces $g_k\subset S^kT^*\ot N$, $k\ge0$, with
$g_0=N$ and $g_k\subset g_{k-1}^{(1)}$, is called a symbolic
system.
 \end{dfn}

If a system of PDEs $\E$ is given as $F_1=0,\dots,F_r=0$, where
$F_i$ are scalar PDEs on $M$, then $T=TM,N\simeq\R^m$ and the
system $g\subset ST^*\ot N$ is given as $f_1=0,\dots,f_r=0$, where
$f_i=\sigma(F_i)$ are symbols of the differential operators at the
considered point (or jet for non-linear PDEs).

With every such a system we associate its Spencer $\d$-complex of
order $k$:
 \begin{multline*}
0\to g_k\stackrel{\d}\longrightarrow g_{k-1}\ot
T^*\stackrel{\d}\longrightarrow
g_{k-2}\ot\La^2T^*\to\dots\\
\to g_i\ot\La^{k-i}T^*\stackrel{\d}\longrightarrow
\dots\stackrel{\d}\longrightarrow g_{k-n}\ot\La^nT^*\to0.
 \end{multline*}

 \begin{dfn}
The cohomology group at the term $g_i\ot\La^jT^*$ is denoted by
$H^{i,j}(g)$ and is called the Spencer $\d$-cohomology of $g$.
 \end{dfn}

Note that $g_k=S^kT^*\ot N$ for $0\le k<r$ and the first number
$r=r_\text{min}(g)$, where the equality is violated is called the
minimal order of the system. Actually the system has several
orders:
 $$
\op{ord}(g)=\{k\in\Z_+\,|\,g_k\ne g_{k-1}^{(1)}\}.
 $$

Multiplicity of an order $r$ is:
 $$
m(r)=\dim g_{r-1}^{(1)}/g_r=\dim H^{r-1,1}(g).
 $$

Hilbert basis theorem implies finiteness of the set of orders
(counted with multiplicities):
 $$
\op{codim}(g):=\dim H^{*,1}(g)=\sum m(r)<\infty.
 $$

Starting from the maximal order of the system $k=r_\text{max}$ we
have: $g_{k+l}=g_k^{(l)}$.

If we dualize the above construction over $\R$, then Spencer
$\d$-differential transforms to a homomorphism over the algebra of
polynomials $ST$ and $g^*=\oplus_i g_i^*$ becomes an $ST$-module.
This module is called a {\em symbolic module\/} and it plays an
important role in understanding PDEs.

In particular, {\em characteristic variety\/}
$\op{Char}^\C(g)\subset \mathbb{P}^\C T^*$ is defined as the
support of this module $\op{Supp}(g^*)=\{[p]:(g^*)_p\ne0\}$ and
the {\em characteristic sheaf\/} $\mathcal{K}$ over it is the
family of vector spaces, which at the point $p\in\op{Char}^\C(g)$
equals the value of the module at this point
$\mathcal{K}_p=g^*/p\cdot g^*$. For more geometric description see
\cite{S,KLV,KL$_2$}.

\section{Commutative algebra approach}

We will study only local solutions of a system of PDEs $\E$, which
we consider in such a neighborhood that type of the symbolic
system does not change from point to point (on equation) in the
sense that dimensions of $g_k$, of the characteristic variety
$\op{Char}^\C(g)$ and of the fibers of $\mathcal{K}$ are the same.

It should be noted that if a system $\E$ is not formally
integrable and $\E'$ is obtained from it by the
prolongation-projection method \cite{K,M2,KL$_2$}, then the
numbers $p,\sigma$ change in this process, i.e. either the
functional dimension or the functional rank decrease. Thus from
now on we suppose the system $\E$ is formally integrable.

The numbers $p,\sigma$ can be described using the methods of
commutative algebra. Recall (\cite{AM}) that by Hilbert-Serre
theorem the sum
 $$
f(k)=\sum_{i\le k}\dim g_i^*
 $$
behaves as a polynomial in $k$ for sufficiently large $k$. This
polynomial is called the {\em Hilbert polynomial\/} of the
symbolic module $g^*$ corresponding to $\E$ and we denote it by
$P_\E(z)$. If $p=\deg P_\E(z)$ and $\sigma=P^{(p)}_\E(z)$, then
the highest term of this polynomial is
 $$
P_\E(z)=\sigma z^p+\dots
 $$
(see \cite{H} for the related statements in algebraic geometry,
the interpretation for PDEs is straightforward).

A powerful method to calculate the Hilbert polynomial is
resolution of a module. In our case a resolution of the symbolic
module $g^*$ exists and it can be expressed via the Spencer
$\d$-cohomology. Indeed, the Spencer cohomology of the symbolic
system $g$ is $\R$-dual to the Koszul homology of the module $g^*$
and for algebraic situation this resolution was found in
\cite{Gr}. It has the form:
 \begin{multline*}
0\to\oplus_q H^{q-n,n}(g)\otimes S^{[-q]}
\stackrel{\varphi_n}\longrightarrow \oplus_q
H^{q-n+1,n-1}(g)\otimes
S^{[-q]} \stackrel{\varphi_{n-1}}\longrightarrow\dots\\
\to\oplus_q H^{q-1,1}(g)\otimes S^{[-q]}
\stackrel{\varphi_1}\longrightarrow \oplus_q H^{q,0}(g)\otimes
S^{[-q]} \stackrel{\varphi_0}\longrightarrow g^*\to0,
 \end{multline*}
where $S^{[-q]}$ is the polynomial algebra on $T^*_xM$ with
grading shifted by $q$, i.e. $S^{[-q]}_i=S^{i-q}T_xM$, and the
maps $\varphi_j$ have degree 0.

Thus denoting $h^{i,j}=\dim H^{i,j}(g)$ and $\tau_\alpha=\dim
S^\alpha TM=\binom{\alpha+n-1}\alpha$ we have:
 $$
\dim g_i=\sum_q \bigl( h^{q,0}\tau_{i-q}-h^{q,1}\tau_{i-q-1}+
h^{q,2}\tau_{i-q-2}-\dots+(-1)^n h^{q,n}\tau_{i-q-n}\bigr).
 $$
Let also $j_\beta=\sum_{\alpha\le\beta}\tau_\alpha=\dim
J^\beta_vM= \binom{\beta+n}n$ be the dimension of the fiber of the
vertical jets $J^\beta_vM$, i.e. the fiber of the jet space
$J^\beta M$ over $M$. Thus we calculate
 $$
\sum_{i\le k}\dim g_i=\sum_q \bigl(
h^{q,0}j_{k-q}-h^{q,1}j_{k-q-1}+ h^{q,2}j_{k-q-2}-\dots\pm
h^{q,n}j_{k-q-n}\bigr).
 $$

Finally we deduce the formula for Hilbert polynomial of the
symbolic module $g^*$
 \begin{multline*}
P_\E(z)=\sum_q \Bigl(
h^{q,0}\tbinom{z-q+n}n-h^{q,1}\tbinom{z-q+n-1}n+\\
+h^{q,2}\tbinom{z-q+n-2}n-\dots +(-1)^n
h^{q,n}\tbinom{z-q}n\Bigr).
 \end{multline*}
Here
 $$
\tbinom{z+k}k=\frac1{k!}(z+1)\cdot(z+2)\cdots(z+k).
 $$
Denote
$S_j(k_1,\dots,k_n)=\sum\limits_{i_1<\cdots<i_j}k_{i_1}\cdots
k_{i_j}$ the $j$-th symmetric polynomial and let also
 $$
s_i^n=\frac{(n-i)!}{n!}\,S_i(1,\dots,n)
 $$
Thus
 $$
s_0^n=1,\quad s_1^n=\frac{n+1}2,\quad
s_2^n=\frac{(n+1)(3n+2)}{4\cdot3!},\quad
s_3^n=\frac{n(n+1)^2}{2\cdot4!},
 $$
 $$
s_4^n=\frac{(n+1)(15n^3+15n^2-10n-8)}{48\cdot5!}\quad \text{ etc.}
 $$

If we decompose
 $$
\tbinom{z+n}n=\sum_{i=0}^n s_i^n\frac{z^{n-i}}{(n-i)!},
 $$
then we get the expression for the Hilbert polynomial
 $$
P_\E(z)=\sum_{i,j,q}(-1)^i h^{q,i} s_j^n
\frac{(z-q-i)^{n-j}}{(n-j)!}\\
=\sum_{k=0}^n b_k\frac{z^{n-k}}{(n-k)!},
 $$
where
 $$
b_k=\sum_{j=0}^k\sum_{q,i}(-1)^{i+j+k}h^{q,i}s_j^n
\frac{(q+i)^{k-j}}{(k-j)!}.
 $$

\section{Calculations for the Solutions space}\label{S6}

We are going to compute the dimensional characteristics of two
important classes of PDEs.

{\bf Involutive systems.} These are such symbolic systems
$g=\{g_k\}$ that all subspaces $g_k$ are involutive in the sense
of Cartan \cite{C,BCG$^3$} (this definition for the symbolic
systems of different orders was introduced in \cite{KL$_5$}).
Thanks to Serre's contribution \cite{GS} we can reformulate this
via Spencer cohomology as follows.

Denote by $g^{|k\rangle}$ the symbolic system generated by all
differential corollaries of the system deduced from the order $k$:
 $$
g^{|k\rangle}_i=\left\{\begin{array}{ll}S^iT^*\ot N,&\text{ for }i<k;\\
g_k^{(i-k)},&\text{ for }i\ge k.\end{array}\right.
 $$
Then the system $g$ is involutive iff $H^{i,j}(g^{|k\rangle})=0$
for all $i\ge k$ (this condition has to be checked for
$k\in\op{ord}(g)$ only), see \cite{KL$_5$}.

In particular, $H^{i,j}(g)=0$ for $i\notin\op{ord}(g)-1$,
$(i,j)\ne(0,0)$, and the resolution for the symbolic module $g^*$
as well as the formula for the Hilbert polynomial of $\E$ become
easier.

Let us restrict for simplicity to the case of systems of PDEs $\E$
of pure first order. Then
 \begin{multline*}
P_\E(z)=h^{0,0}\tbinom{z+n}n-h^{0,1}\tbinom{z+n+1}{n+1}+
h^{0,2}\tbinom{z+n+2}{n+2}-\dots\\ = b_1\frac{z^{n-1}}{(n-1)!}+
b_2\frac{z^{n-2}}{(n-2)!}+\dots+b_0.
 \end{multline*}
Vanishing of the first coefficient $b_0=0$ is equivalent to
vanishing of Euler characteristic for the Spencer $\d$-complex,
$\chi=\sum_i(-1)^ih^{0,i}=0$, and this is equivalent to the claim
that not all the covectors from $^\C T^*\setminus0$ are
characteristic for the system $g$.

The other numbers $b_i$ are given by the general formulas from the
previous section, but they simplify in our case. For instance
 $$
b_1=\tfrac{n+1}2\,b_0-\sum(-1)^ih^{0,i}i=\sum(-1)^{i+1}i\cdot
h^{0,i}.
 $$

If $\op{codim}\op{Char}^\C(\E)=n-p>1$, then $b_1=0$ and in fact
then $b_i=0$ for $i<n-p$, but $b_{n-p}=\z$.

 \begin{theorema}
If $\op{codim}\op{Char}^\C(\E)=n-p$, then the functional rank of
the system equals
 $$
\sigma=\sum_i (-1)^ih^{0,i}\frac{(-i)^{n-p}}{(n-p)!}.
 $$
 \end{theorema}

 \begin{proof}
Indeed one successively calculate the coefficients using the
formula
  $$
b_k=\sum_i\sum_{\alpha=0}^k (-1)^{i+\alpha}h^{0,i}s_{k-\alpha}^n
\frac{i^\alpha}{\alpha!}
 $$
and notes that $b_k$ equals to the displayed expression plus a
linear combination of $b_{k-1},\dots,b_0$. The claim follows.
 \end{proof}

One can extend the above formula for general involutive system and
thus compute the functional dimension and functional rank of the
solutions space (some interesting calculations can be found in
classical works \cite{J,C}).

{\bf Cohen-Macaulay systems.} A symbolic system $g$ (and the
respective PDEs system $\E$) is called Cohen-Macaulay
(\cite{KL$_2$}) if the corresponding symbolic module $g^*$ is
Cohen-Macaulay, i.e. (see \cite{M1,E} for details)
 $$
\dim g^*=\op{depth}g^*.
 $$

Consider an important partial case (we formulate the definition
only for symbolic systems; PDEs are treated in \cite{KL$_4$}):

 \begin{dfn}
A symbolic system $g\subset ST^*\ot N$ ($n=\dim T$, $m=\dim N$) of
$\op{codim}(g)=r$ is called a {\em generalized complete
intersection\/} if
 \begin{itemize}
\item $m\le r<n+m$;
\item $\op{codim}_\C\op{Char}^\C(g)=r-m+1$;
\item $\dim\mathcal{K}_x=1$
 $\forall x\in\op{Char}^\C(g)\subset P^\C T^*$.
 \end{itemize}
 \end{dfn}

Formal integrability of such systems are given by the
compatibility conditions expressed via brackets (for scalar
systems \cite{KL$_1$,KL$_3$}) or multi-brackets (for vector
systems \cite{KL$_4$}). In this case we can calculate Cartan genre
and integer directly.

 \begin{theorema}
Let $\E$ be a system of generalized complete intersection type and
suppose it is formally integrable. Then the functional dimension
of $\op{Sol}(\E)$ is
 $$
p=m+n-r-1
 $$
and the functional rank is
 $$
\sigma=S_{r-m+1}(k_1,\dots,k_r)=\sum\limits_{i_1<\dots<i_{r-m+1}}
k_{i_1}\cdots k_{i_{r-m+1}},
 $$
the $l$-th symmetric polynomial of the orders $k_1,\dots,k_r$ of
the system.
 \end{theorema}

Note that if the last requirement in the definition of generalized
complete intersection is changed to $\dim\mathcal{K}_x=d$
everywhere on the characteristic variety, then the functional rank
will be multiplied by $d$:
 $$
\sigma=d\cdot S_{r-m+1}(k_1,\dots,k_r).
 $$
However the formal integrability criterion for generalized
complete intersections is proved in \cite{KL$_4$} under assumption
that $d=1$.

 \begin{proof}
We shall consider the case of a system $g$ of a pure order:
$k_1=\dots=k_r=k$, $k_i\in\op{ord}(g)$. The case of different
orders is similar and will appear elsewhere.

The formula for functional dimension $p$\, follows directly from
the definition of generalized complete intersection. Let's
calculate $\sigma$.

We can use interpretation of the Cartan integer $\sigma$ from
\S\ref{Cn}. Recall that characteristic variety $\op{Char}^\C(g)$
is the locus of the characteristic ideal $I(g)=\op{Ann}(g)$, which
the the annihilator of $g^*$ in $ST$.

Since the module is represented by the matrix with polynomial
entries (each differential operator $\Delta_i$ giving a PDEs
system $\E$ is a column $\Delta_{ij}$, $1\le i\le r$, $1\le j\le
m$; so that their union is a $m\times r$ matrix $M(\Delta)$), its
annihilator is given by the zero Fitting ideal (in fact, here we
use the condition on grade of the ideal:
$\op{depth}\op{Ann}(g)=r-m+1$, which follows from the conditions
of the above definition).

This ideal $\op{Fitt}_0(g)$ is generated by all determinants of
$m\times m$ minors of the corresponding to $M(\Delta)$ matrix of
symbols $M(\sigma_\Delta)$. These minors are determined by a
choice of $m$ from $r$ columns, so that there are $\tbinom{r}m$
determinants and each is a polynomial of degree $k^{r-m+1}$.

However not all the minors are required to determine
$\op{Char}^\C(g)$ and this is manifested by the fact, that we sum
$\tbinom{r}{m-1}$ degrees $k^{r-m+1}$ to get the functional rank
$\sigma$. The easiest way to explain this is via the Hilbert
polynomial of the symbolic module $g^*$.

This can be calculated since under the assumption of generalized
complete intersection $g^*$ possesses a resolution in the form of
Buchsbaum-Rim complex (see \cite{KL$_4$}):
 \begin{multline*}
0\to S^{r-m-1}V^\star\ot\Lambda^r U\to
S^{r-m-2}V^\star\otimes\Lambda^{r-1}U \to\\
\dots\to\Lambda^{m+1}U\to U\to V\to g^*\to0,
 \end{multline*}
where $V\simeq ST\otimes N^*$ (recall that $\dim N=m$ and
$g\subset ST^*\otimes N$) and $U=\underbrace{ST\oplus\dots\oplus
ST}_{\mbox{$r$ terms}}$. Star $\star$ means dualization over $ST$
and the tensor products are over $ST$ as well.

Now the claim follows from the detailed investigation of degrees
of the homomorphisms in the above exact sequence. To see this we
suppose at first that $r=m+n-1$ and use the following assertion.

 \begin{lem}
The following combinatorial formula holds:
 \begin{multline*}
m\tbinom{n+k(n+m-1)}{n}-(n+m-1)\tbinom{n+k(n+m-2)}n\\
+\sum_{j=1}^{n-1}(-1)^{j-1}\tbinom{j+m-2}{m-1}
\tbinom{n+m-1}{j+m}\tbinom{(k+1)n-k(1+j)}n =\tbinom{n+m-1}n k^n.
 \end{multline*}
 \end{lem}

We would like to comment and interpret the sum on the left hand
side of this formula. In our case the system is of finite type
($g^*$ has finite dimension as a vector space) and
$\sigma=\sum\dim g_i$ (the sum is finite).

Stabilization of the symbol occurs at the order $i=\sum
k_i-1=k(n+m-1)-1$: $g_i=0$. So we prolong $\E$ to the jets of
order $k(n+m-1)$ and the first term is just $\dim
J^{k(n+m-1)}_v(M,N)$.

The next term is due to the fact that $\E\subset J^k(M,N)$ is
proper. It is given by $r=n+m-1$ equations of order $k$, we which
we differentiate up to $k(n+m-2)$ times along all coordinate
directions (prolongation).

There are relations between these derivatives. These are
compatibility conditions (1-syzygy of $g^*$), which appear in the
form of multi-brackets \cite{KL$_4$}, in our case this bracket
uses $(m+1)$-tuples of $\Delta_i$.

There are in turn relations among relations (2-syzygy of $g^*$),
which are identities between multi-brackets (these we call
generalized Pl\"ucker identities, to appear soon), in our case
these latter use $(m+2)$-tuples of the defining operators
$\Delta_i$ etc.

Due to exact form of the relations (higher syzygies) we get
factors $\tbinom{j+m-2}{m-1}$ in the summations formula of the
lemma.

In the case $r<n+m-1$ we should perform a reduction, which is
possible by Theorem A \cite{KL$_2$}. Then the functional dimension
$p$ grows, but the functional rank remains the same and the
previous calculation works.
 \end{proof}

\section{Examples}

Here we show some examples demonstrating the above results.

{\bf 1.} {\em Intermediate integral} of a system $\E\subset
J^k\pi$ is such a system $\tilde\E\subset J^{\tilde k}\pi$ that
$\tilde k<k$ and $\E\subset\tilde\E{}^{(k-\tilde k)}$ (where
$\E^{(i)}$ is the $i$-th prolongation of the system). Since every
solution to the system $\E$ is a solution to
$\tilde\E{}^{(k-\tilde k)}$ we conclude: Whenever the functional
dimension $p>0$, we have $\tilde p=p$ and $\tilde\sigma=\sigma$.

Indeed the solutions of $\tilde\E{}^{(k-\tilde k)}$ form a
finite-dimensional parametric family, such that solutions of
$\tilde\E$ appear for some fixed values of parameters (because we
differentiate with respect to all variables to obtain the
prolongation). Thus the number of functions of $p>0$ variables, on
which a general solution depends, will not be altered.

{\bf 2.} If the PDEs system $\E$ is underdetermined, then $p=n$
and $\sigma\ge1$. Indeed, $\sigma$ is precisely the
under-determinacy degree, i.e. the minimal number of unknown
functions that should be arbitrarily fixed to get a determined
system. We assume we can do it to get a formally integrable
system. If underdetermined system is not formally integrable,
compatibility conditions can turn it into determined or
over-determined and then decrease $p$ and change $\sigma$.

A nice illustration is the Hilbert-Cartan system
 $$
z'(x)=(y''(x))^2.
 $$
It has $p=1$, $\sigma=1$. But even though a general solution
depends on one function of one variable, it cannot be represented
in terms of a function and its derivatives only (Hilbert's
theorem).

{\bf 3.} As we noticed earlier the similar situation happens to
overdetermined system: If $\E$ is not formally integrable, and
$\tilde\E$ is obtained from $\E$ by prolongation-projection
technique (sometimes it is said that $\tilde\E$ is the involutive
form of $\E$, but this is not true, only a certain prolongation of
$\tilde\E$ is), then $\tilde p<p$ or [$\tilde p=p$ and
$\tilde\sigma<\sigma$]. Indeed, supplement of additional equations
shrinks the solution space.

For instance if we consider two second-order scalar differential
equations on the plane
 \begin{gather*}
F\bigl(x,y,u(x,y),Du(x,y),D^2u(x,y)\bigr)=0,\\
G\bigl(x,y,u(x,y),Du(x,y),D^2u(x,y)\bigr)=0,
 \end{gather*}
such that $F$ and $G$ have no common complex characteristics, then
the compatibility condition of this system $\E$ can be expressed
via the Mayer bracket (\cite{KL$_1$}): $H=[F,G]_\E$. If $H=0$,
then $p=0$, $\sigma=4$. If $H\ne0$, then $p=0$ and $\sigma\le3$,
the equality being given by the Frobenius condition for the system
$\tilde\E=\{F=0,G=0,H=0\}$.

If the system has one common characteristic and is compatible, we
have: $p=1$, $\sigma=1$. Pairs of such systems are basic examples
of Darboux integrability.

{\bf 4.} Evolutionary equations $u_t=L[u]$ provide interesting
examples, which usually "{}contradict"\ the theory. Consider for
instance the heat equation
 $$
u_t=u_{xx}.
 $$
It is formally integrable and analytic. We can try to specify the
initial condition $u|_{t=0}=\varphi(x)$ and then solve the Cauchy
problem, so that we get $p=1$, $\sigma=1$. On the other hand we
can let $u|_{x=0}=\psi_0(t),u_x|_{x=0}=\psi_1(t)$ and then get
$p=1$, $\sigma=2$.

If we calculate the numbers using our definitions of functional
dimension and functional rank (for instance, via Hilbert
polynomial), it turns out that the second approach is correct.
Indeed with the first idea we come into trouble with certain
Cauchy data: Let, for instance, $\varphi(x)=(1-x)^{-1}$, which is
an analytic function around the origin. Then the analytic solution
should have the series
 $$
u(t,x)\doteqdot\frac1{1-x}+\frac21\frac{t}{(1-x)^3}
+\frac{4!}{2!}\frac{t^2}{(1-x)^5}+\dots
+\frac{(2n)!}{n!}\frac{t^n}{(1-x)^{2n+1}}+\dots
 $$
which diverges everywhere outside $t=0$. The reason why the second
approach provides no problem is because the line $\{x=0\}$ is
non-characteristic and we can solve our first order PDE by the
classical method of Cauchy characteristics.

Remark however that in the standard courses of mathematical
physics the heat equation is solved with the first approach (by
Fourier method). How is it possible?

Explanation is that we solve the heat equation then only for
positive time $t\ge0$. Doing the same method in negative direction
blows up the solutions immediately (heat goes rapidly to
equilibrium, but we cannot predict even closest past)! We here are
interested in the solutions, which exist in an open neighborhood
of the origin (like in Cauchy-Kovalevskaya theorem), and this
contradicts the first approach.

{\bf 5.} Similar problems arise with Cauchy problems in other PDEs
systems: one usually applies reduction or fixes gauge, but this
can change dimensional characteristics.

For instance, consider the Cauchy problem for the Einstein vacuum
equations, which is a system of 10 PDEs of 10 unknown functions.
The system is over-underdetermined (i.e. it has compatibility
conditions). In wave gauge \cite{CB} its solution depends on
several functions on a 3-dimensional space, which are subject to
constraint equations, so that $p=2$. On the other hand, the
original Einstein system is invariants under diffeomorphisms and
this yields $p=4$.

One should also be careful with Cauchy data in higher order, since
then the definition of characteristics becomes more subtle, see
\cite{KL$_5$}.

{\bf 6.} Consider a system $\E$, which describes automorphisms of
a given geometric structure. The corresponding symbolic system is
$g\subset ST^*\otimes T$. The automorphism group has maximal
dimension iff the system is formally integrable. Consider the
examples, when the geometric structure is symplectic, complex or
Riemannian (all these structures are of the first order).

Let at first $g$ be generated by $g_1=\op{sp}(n)\subset T^*\ot T$.
Our tangent space $T=T_xM$ is equipped with a symplectic structure
$\omega$, and we can identify $T^*\stackrel\omega\simeq T$ and we
get $g_1=S^2T^*\subset T^*\otimes T^*$. The prolongations are
$g_i=S^{i+1}T^*\subset S^iT^*\otimes T$.

The system is easily checked to be involutive and the only
non-vanishing Spencer $\d$-cohomology groups are
 $$
H^{0,i}(g)=\La^{i+1}T^*.
 $$
Then one checks that the Euler characteristic is $\chi=1\ne0$ and
so $b_0\ne0$. Thus the functional dimension is $p=n$. Indeed the
characteristic variety is $\mathbb{P}^\mathbb{C}T^*$ because each
non-zero covector $p$ is characteristic: $p^2\in g_1\simeq
S^2T^*$. Next by a theorem from \S\ref{S6} one calculates the
functional rank
 $$
\sigma=\sum_{i=0}^{n-1}(-1)^i\binom{n}{i+1}\frac{(-i)^0}{0!}=\chi=1.
 $$
This result is easy to verify: an infinitesimal symplectic
transformation has a generating function (Hamiltonian) and so it
is determined by one function of $n$ variables.

If we turn to (almost) complex structures $J$ on $M$, then
$g_1=\op{gl}(\frac n2,\mathbb{C})=T^*\otimes_\mathbb{C} T$ (space
of $\mathbb{C}$-linear endomorphisms of $T$) and the prolongations
are $g_i=S^i_\mathbb{C}T^*\otimes_\mathbb{C} T$.

The characteristic variety is proper and one calculates that
$p=\frac n2$, $\sigma=n$. The system is again involutive. The
second Spencer cohomology is
 $$
H^{0,2}(g)=\La^2_\mathbb{C}T^*\otimes_{\bar{\mathbb{C}}}T,
 $$
which is the space of $\mathbb{C}$-antilinear skew-symmetric
$(2,1)$ tensors (Nijenhuis tensors).

The last example is the algebra of Riemannian isometries (i.e. $T$
is equipped with a Riemannian structure) of a Riemannian metric
$q$ on $M$. The symbol is $g_1=o(n)$ and the prolongations are
zero $g_2=g_3=\dots=0$.

This system is not involutive. For instance,
 $$
H^{1,2}(g)=\op{Ker}\bigl(S^2\La^2T^*\to\La^4T^*\bigr)
 $$
(the space of Riemannian curvatures) is non-zero (for $n=\dim
T>1$). Since the system is of finite type, the characteristic
variety is empty and $p=0$. The general solution (isometry)
depends on $\sigma=\frac{(n+1)n}2$ constants.

We recall, that the above dimensional conclusions are correct if
the system $\E$ is integrable, otherwise the space $\op{Sol}(\E)$
shrinks. In the above examples this means: the form $\omega$ is
closed (with just non-degeneracy we have almost-symplectic
manifold); the structure $J$ is integrable (Nijenhuis tensor $N_J$
vanishes); the manifold $(M,q)$ has constant sectional curvature
(so it is a spacial form).

\end{document}